\documentclass[conference]{IEEEtran}
\usepackage[backend=bibtex,bibstyle=ieee,citestyle=numeric-comp]{biblatex}
\bibliography{references}
\IEEEoverridecommandlockouts
\usepackage{amsmath,amssymb,amsfonts}
\usepackage{algorithmic}
\usepackage{graphicx}
\usepackage{graphics} 
\usepackage{epsfig} 
\usepackage{mathptmx} 
\usepackage{times} 
\usepackage{amssymb}  
\usepackage{amsmath} 
\usepackage{newtxtext, newtxmath}
\usepackage{todonotes}
\usepackage{booktabs}
\usepackage{lipsum}
\usepackage{hyperref}
\usepackage[font=small]{caption}
\usepackage{subcaption}

\DeclareMathOperator*{\argmin}{arg\,min}
\hypersetup{bookmarksopen,bookmarksnumbered,
pdfpagemode=UseOutlines,
colorlinks=true,
urlcolor=blue,
linkcolor=blue,
anchorcolor=blue,
citecolor=blue,
filecolor=blue,
menucolor=blue
}
\usepackage{flushend}
\usepackage{cleveref}
\usepackage{textcomp}
\usepackage{xcolor}
\def\BibTeX{{\rm B\kern-.05em{\sc i\kern-.025em b}\kern-.08em
    T\kern-.1667em\lower.7ex\hbox{E}\kern-.125emX}}
\begin{document}

\title{Mitigating Landside Congestion at Airports through \\ Predictive Control of Diversionary Messages\\
\thanks{*Pacific Northwest National Laboratory is operated by Battelle Memorial Institute for the U.S. Department of Energy under Contract No. DE-AC05-76RL01830. This work was supported by the U.S. Department of Energy Vehicle Technologies Office.}
}
\makeatletter
\newcommand{\linebreakand}{%
  \end{@IEEEauthorhalign}
  \hfill\mbox{}\par
  \mbox{}\hfill\begin{@IEEEauthorhalign}
}
\makeatother

\author{\IEEEauthorblockN{Nawaf Nazir}
\IEEEauthorblockA{\textit{Pacific Northwest National Lab} \\
Richland, USA \\
nawaf.nazir@pnnl.gov}
\and
\IEEEauthorblockN{Soumya Vasisht}
\IEEEauthorblockA{\textit{Pacific Northwest National Lab} \\
Richland, USA \\
soumya.vasisht@pnnl.gov}
\and
\IEEEauthorblockN{Shushman Choudhury}
\IEEEauthorblockA{\textit{Lacuna Technologies} \\
Palo Alto, USA \\
shushman.choudhury@lacuna.ai}
\linebreakand
\IEEEauthorblockN{Stephen Zoepf}
\IEEEauthorblockA{\textit{Lacuna Technologies} \\
Palo Alto, USA \\
stephen.zoepf@lacuna.ai}
\and
\IEEEauthorblockN{Chase Dowling}
\IEEEauthorblockA{\textit{Pacific Northwest National Lab} \\
Richland, USA \\
chase.dowling@pnnl.gov}
}

\maketitle

\begin{abstract}
We present a data-driven control framework for adaptively managing landside congestion at airports. Ground traffic significantly impacts airport operations and critical efficiency, environmental, and safety metrics. Our framework models a real-world traffic intervention currently deployed at Seattle-Tacoma International Airport (SEA), where a digital signboard recommends drivers to divert to Departures or Arrivals depending on current congestion. We use measured vehicle flow/speed and passenger volume data, as well as time-stamped records of diversionary messages, to build a macroscopic system dynamics model. We then design a model predictive controller that uses our estimated dynamics to recommend diversions that optimize for overall congestion. Finally, we evaluate our approach on $50$ real-world historical scenarios at SEA where no diversions were deployed despite significant congestion. Our results suggest that our framework would have improved speed in the congested roadway by up to three times and saved between $20$ and $80$ vehicle-hours of cumulative travel time in every hour of deployment. Overall, our work emphasizes the opportunity of algorithmic decision-making to augment operator judgment and intuition and yield better real-world outcomes.
\end{abstract}

\begin{IEEEkeywords}
Airports, Traffic congestion, Data-driven control, Predictive control, Diversionary messages
\end{IEEEkeywords}

\section{Introduction}
\label{sec:intro}


The problem of ground traffic congestion at airports has persisted for decades and major airports have invested considerable time and money to alleviate it~\cite{brink1975identification}. The policies enacted range from hiring ground staff that resolve specific incidents to multi-year planning and infrastructure projects~\cite{shriner1999evaluating}. Targeted and effective strategies for managing landside congestion can greatly influence overall quality of service by reducing travel times, idle vehicle emissions, and adverse events, and improving traveler satisfaction~\cite{muller1991framework,failla2014exploring}.

Our paper focuses on a specific adaptive intervention: messages on a digital signboard that recommend incoming vehicles to divert from Departures to Arrivals when the former is congested (and vice versa). We do so because such signboard systems are already in use in certain airports, where the diversions are deployed based on operator intuition and judgment rather than algorithmic decision-making (e.g., Seattle-Tacoma International Airport, or SEA, whose real-world data we use to build a macroscopic traffic model and to evaluate our approach). But our methodology could apply to other adaptive congestion management tactics at airports, and even more widely to other ground traffic problems for which such data are available.

The most relevant prior research on this problem used so-called microscopic, i.e., high-fidelity simulations that modeled the behaviour of individual vehicles~\cite{bender1997simulating,duncan2000development,ugirumurera2021modeling}. Such approaches can represent detailed scenarios and evaluate them comprehensively, but are typically data and time intensive to create and computationally expensive to query. They are far more suitable for assessing the long-term effects of operational changes than for producing real-time interventions. A mesoscopic airport curbside model addressed some but not all of the challenges of microscopic simulations, and focused on analysing various curbside policies rather than optimizing them~\cite{harris2017mesoscopic}. In contrast, to make real-time interventions feasible, we use a macroscopic model with only aggregate traffic information.

\begin{figure}[t]
    \centering
    \includegraphics[width=0.9\columnwidth]{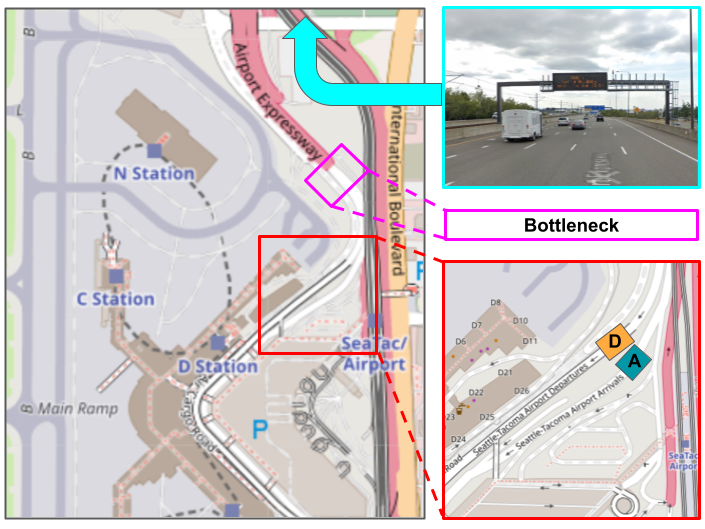}
    \caption{An overhead view of the physical layout of Seattle-Tacoma International Airport, whose real-world data we use. We illustrate the variable messaging sign, traffic bottlenecks, Arrival and Departure drives, and the locations of respective speed/flow sensors, labeled A and D. The figure is used from our previous work~\cite{vasisht2022Estimating}.}
    \label{fig:seatac_map}
\end{figure}

Our key idea is to formulate the problem of adaptive congestion management using diversionary messages as one of model predictive control or MPC, a widely successful framework for adaptively controlling complex systems~\cite{camacho2013model}. Using real-world historical data on vehicle flow and speed and the displayed messages, we estimate a system dynamics model. This model captures both the natural evolution of the macroscopic traffic state as well as its response to interventions. We then design a discrete-time MPC formulation where the control variables are the two kinds of interventions (Departures to Arrivals, and vice versa) and the predictive model is our estimated system dynamics.

\begin{figure*}[t]
    \centering
    \includegraphics[width=0.8\textwidth]{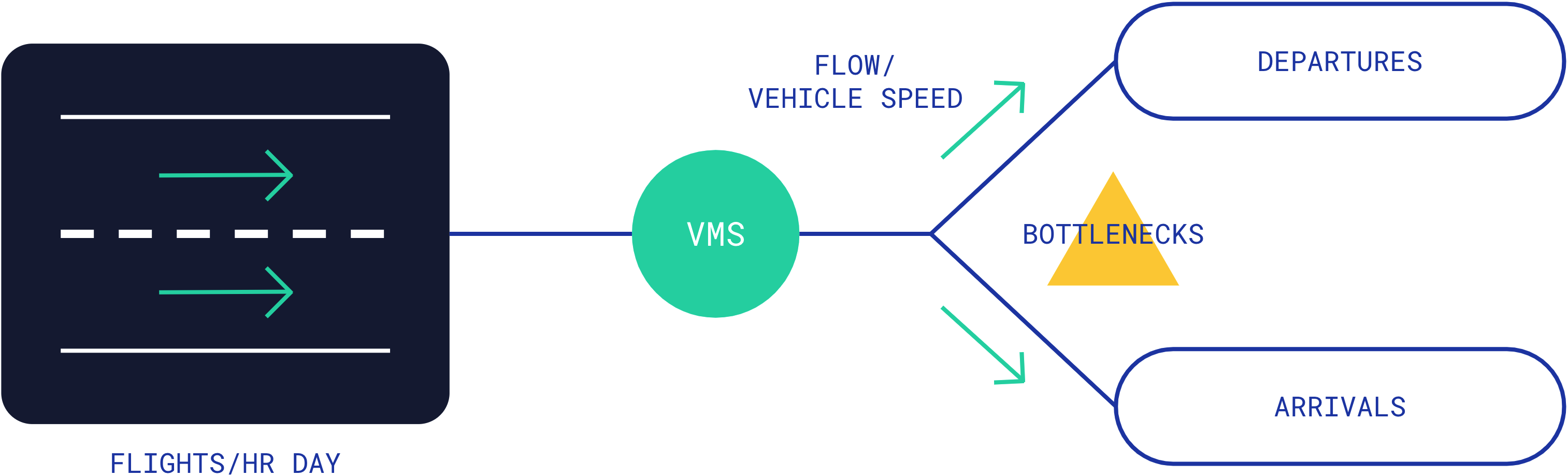}
    \label{fig:seatac-illustr}
\end{figure*}

To evaluate our methodology, we use $50$ historical scenarios from Seattle-Tacoma airport, where one of Departures/Arrivals was congested but no diversionary messages were deployed. Our simulations suggest that intervening adaptively with our framework would have mitigated congestion in all cases. The speed in the congested roadway would have improved significantly (by up to three times), and between $20$ to $80$ vehicle-hours of cumulative travel time would have been saved in every hour of intervention.

\section{Related Work}
\label{sec:background}


This section reviews three areas of work that we build upon. The first is the general idea of using model predictive control and related ideas from optimal control theory to manage ground traffic. The second and third are related but distinct approaches to understanding and improving airport curbside management: modeling and analysis of existing operations and detailed simulations of various potential long-term policies to guide future decision-making.

\subsection{Model Predictive Control for Traffic Management}

The idea of using MPC to solve variants of ground traffic management has been studied in the literature~\cite{gartner1984development}. A comprehensive summary of early approaches can be found in the PhD thesis~\cite{hegyi2004model}. The majority of existing work differs from ours in the specific real-world control strategy that they encode: the timing cycle of traffic signals~\cite{ye2019survey}; so-called ramp-metering to limit inflow of vehicles into a freeway at rush hour~\cite{bellemans2006model}; dynamic curb allocation in urban areas~\cite{Nazir2022curbs}; and perimeter control that applies the general idea of ramp-metering to a large-scale urban network~\cite{sirmatel2017economic,sirmatel2021stabilization}. None of these strategies apply to balancing and diverting inflow between two facilities, each with their unique traffic characteristics. Some recent research exist on the traffic control community considers our flavor of problem, e.g., a theoretical model of flow exchanges in multi-reservoir systems~\cite{mariotte2019flow} and dynamically reorganizing and rearranging inflow on congested lanes and roadways~\cite{anderson2020dynamic,huang2021evaluation}. But they do not consider the complex time-varying nature of airport traffic and do not develop and evaluate their approach on real traffic and intervention data.

\subsection{Understanding Airport Landside Operations}

Since as far back as the late 1960s, travelers to airports have reported facing delays and frustration due to ground traffic congestion and bottlenecks~\cite{brink1975identification}. Over the decades, the literature on evaluating and analysing landside facilities has covered a number of areas: general terminal operations~\cite{tovsic1992review,parr2013evaluating}, access management~\cite{shriner1999evaluating,budd2011airport}, level of service~\cite{muller1991framework,correia2004evaluating}, and customer experience for both quantitative metrics like dwell time~\cite{kim2004estimating} and subjective ones like satisfaction~\cite{failla2014exploring}. Two papers from this line of work are more directly related to our specific problem. One studied how advance signs (albeit for airlines, not diversions) guide and improve vehicle search times~\cite{kichhanagari2002airport} and the other compared the throughput characteristics of pickup to dropoff facilities (of which Arrivals and Departures are canonical examples)~\cite{kleywegt2022throughput}. This body of work is useful for developing context and understanding the general problems and opportunities to improve, but they have no notion of optimization or control that we can use in a computational framework.

\subsection{Landside Management Policies for Airports}

In addition to understanding the problem of landside operations at airports, the research community has been trying to address it as well. Numerical simulations have emerged as a natural and effective tool of choice; an early example is of Las Vegas McCarran International Airport trying to redesign its ground transportation access system~\cite{bender1997simulating}. A related team of authors then used discrete event simulation to model curbside vehicular traffic at airports and applied it to Austin-Bergstrom International Airport~\cite{tunasar1998modeling}. Most early efforts relied on simple bespoke models, e.g., dual roadways for loading/unloading and circulation~\cite{duncan2000development} and a mathematical queueing framework~\cite{fisher2010acrp}; such approaches do not benefit and improve from real-time data about the airports. More recent works have harnessed vastly improved computational resources to develop high fidelity microsimulations that are carefully calibrated to facility data; two notable examples are a study on improving taxi drop-off throughput at a terminal~\cite{yang2020achieving} (they modeled a railway station in Nanjing but the principles hold), and a comprehensive modeling framework of Dallas Fort-Worth International Airport~\cite{ugirumurera2021modeling}. Such detailed simulations can be quite realistic (though they require their fair share of assumptions) and are well-suited for long-term strategic planning. But their data and computation intensive nature is what makes them unsuitable to be used in the loop of a real-time control algorithm.

\section{Methodology}
\label{sec:method}

In this section, we summarize our problem formulation and approach, to emphasize the various interacting components. We operate a variable message signboard (VMS) that incoming drivers observe before they decide whether to go to the Departures or Arrivals roadways (see \Cref{fig:seatac-illustr}). If we infer that either roadway is significantly more congested than the other, we can recommend that drivers heading to the congested roadway divert to the other one. Our problem is to develop a controller that computes a good sequence of diversions for the current and expected traffic conditions. The controller's high-level objective is to mitigate congestion, i.e., improve travel speeds in the congested roadway.

The specifics of our problem formulation depend on the real-world data from SEA that we use to demonstrate and evaluate our approach. We have measurements of vehicle flow (volume per unit time) and average speed over 4 months of the year 2022, binned to 15-minute intervals. The hourly departure and arrival average speeds over a 24-hour horizon are depicted in~\Cref{fig:dep_speed} and~\Cref{fig:arr_speed} respectively. Here, each dot represents the average speed of all vehicles at the terminal within the 15-minute interval. The speed data highlights the difference in congestion times between arrivals and departures and hence the opportunity to divert traffic from one terminal to another when necessary. We also have time-stamped records of messages deployed by SEA, recommending traffic to divert from one terminal to the other. In recent work, we inferred the response rates of incoming traffic to the suggestions~\cite{vasisht2022Estimating}. Here, we will build upon our understanding of the system and design a controller that mitigates congestion.

\begin{figure*}[t]
    \centering
    \begin{subfigure}{0.6\columnwidth}
        \centering
        \includegraphics[width=\columnwidth]{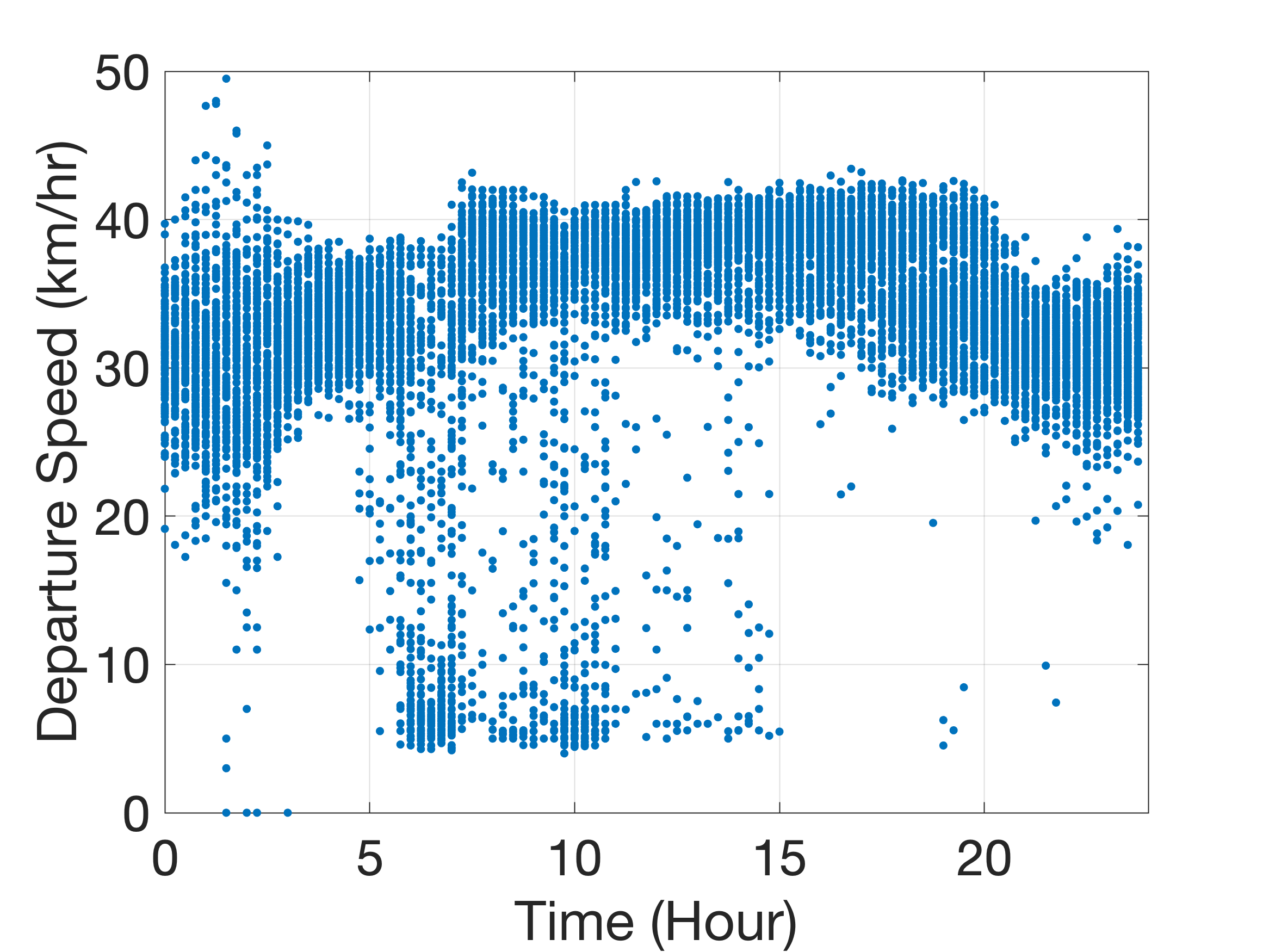}
        \caption{}
        \label{fig:dep_speed}
    \end{subfigure}
    \begin{subfigure}{0.6\columnwidth}
        \centering
        \includegraphics[width=\columnwidth]{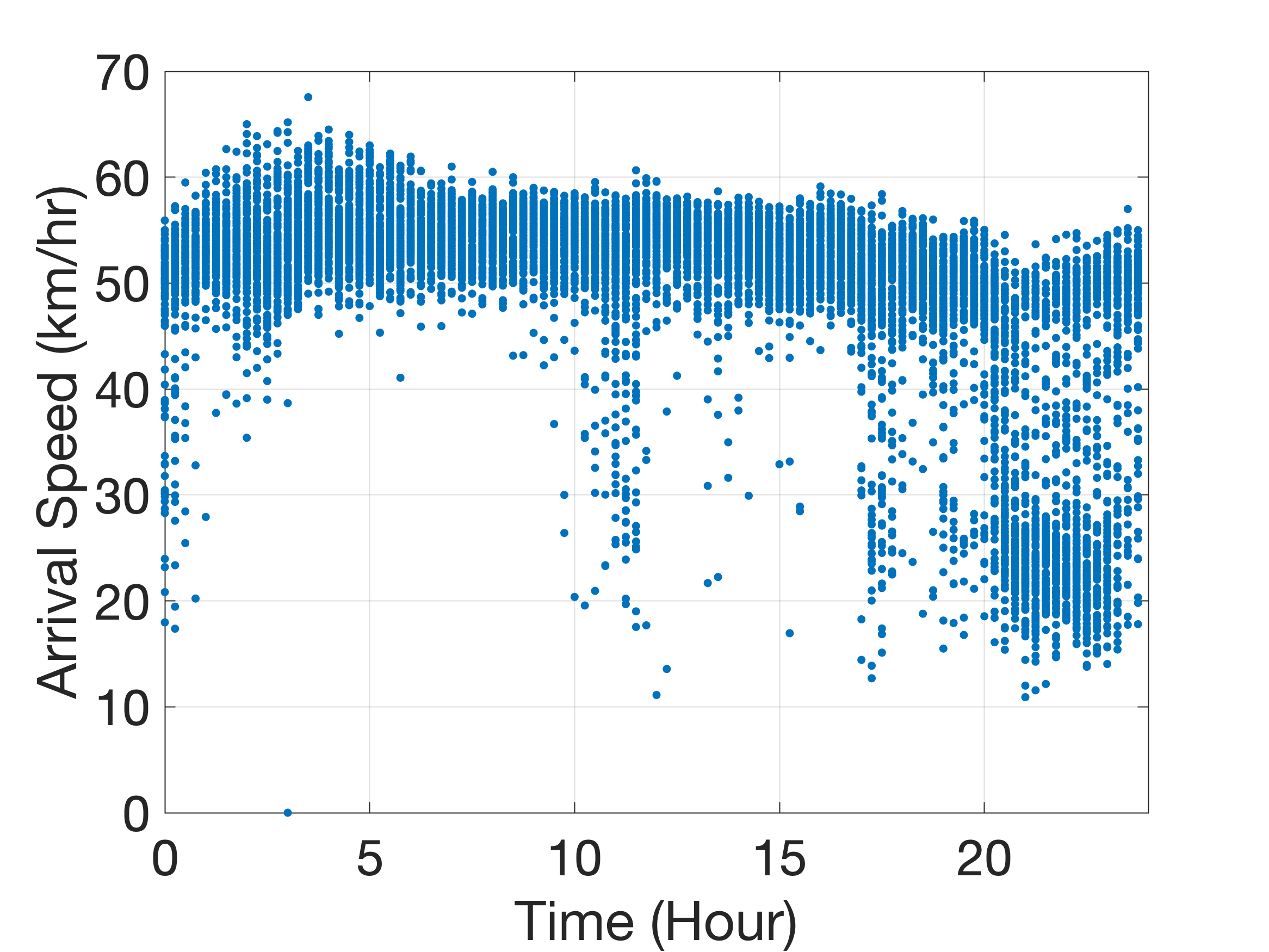}
        \caption{}
        \label{fig:arr_speed}
    \end{subfigure}
    \begin{subfigure}{0.6\columnwidth}
        \centering
        \includegraphics[width=\columnwidth]{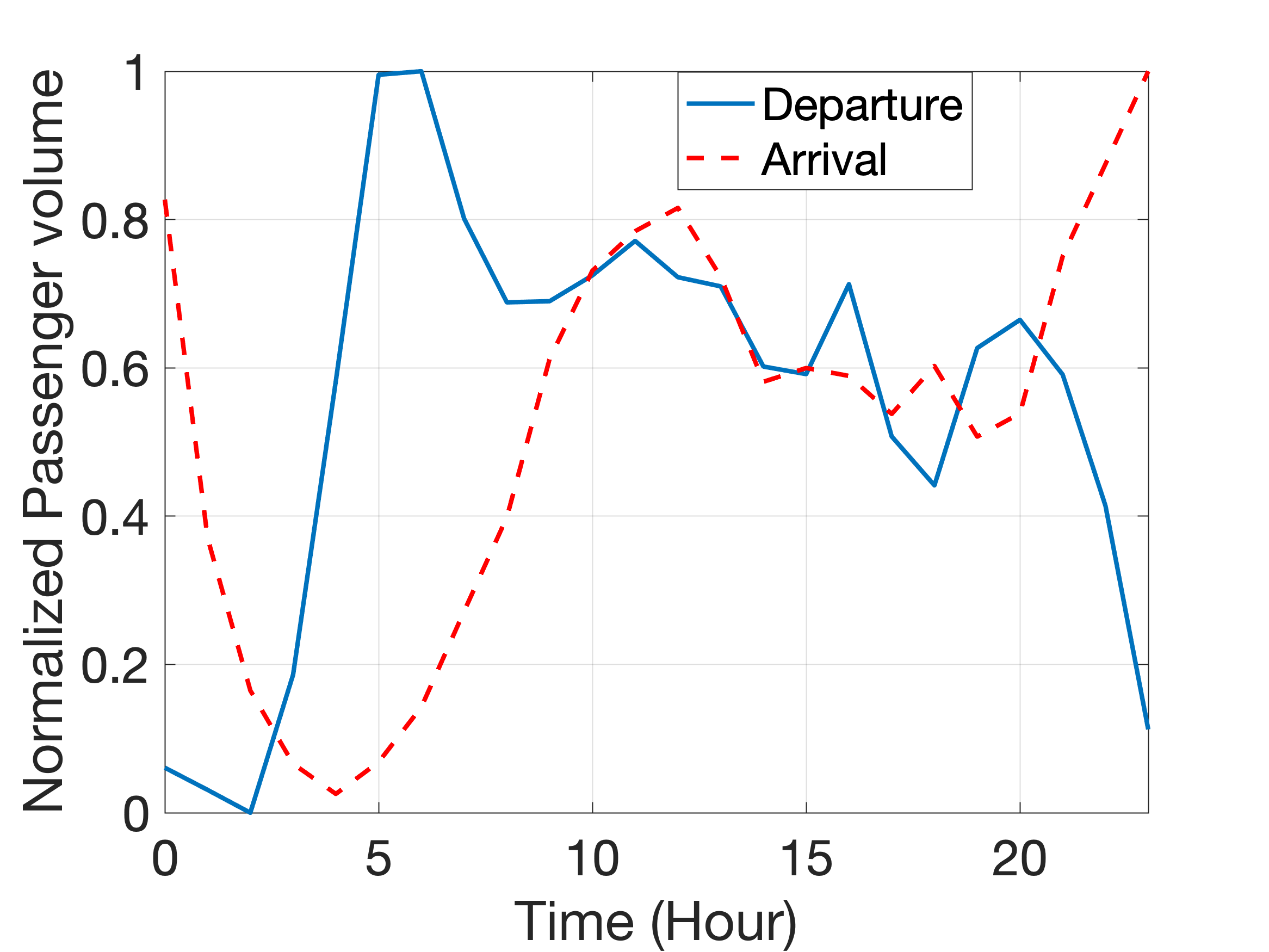}
        \caption{}
        \label{fig:pass-volume-norm}
    \end{subfigure}
    \caption{Average vehicle speeds at (a) Departures and (b) Arrivals and (c) normalized cumulative passenger volumes.}
    \label{fig:arr_dep_speeds}
\end{figure*}
Our approach consists of two key components. First, we estimate a model of the traffic dynamics from the available data, both the natural variation and the response to interventions. We use least-squares regression to fit a linear discrete-time dynamical model. Second, we use the estimated dynamics model in the loop of a model predictive controller. The controller optimizes over sequence of binary control inputs that encode the two kinds of diversions.

\section{Estimating Macroscopic Traffic Dynamics}
\label{sec:dynamics}
Consider a linear discrete-time dynamic model to represent the evolution of airport traffic dynamics:
\begin{align}\label{eq:dynamic_model}
    x_{k+1}=Ax_{k}+Bu_{k}
\end{align}
where $x_k\in \mathbb{R}^{4\times 1}$, $u_k\in \mathbb{R}^{4\times 1}$, $A\in \mathbb{R}^{4\times 4}$ is the auto-regressive co-efficient matrix and $B\in \mathbb{R}^{4\times 4}$ is the control co-efficient matrix. The state $x_{k}$ describes Departures flow and speed (DF/DS) and Arrivals flow and speed (AF/AS) at time $k$ and can be concisely denoted as $x_k=\begin{bmatrix} \textrm{DF, DS, AF, AS} \end{bmatrix}^\top$. The control input $u_k$ represents the treatment applied at time step $k$ for diverting traffic from Departures to Arrivals (TD) or vice versa (TA). We also consider the Departures and Arrivals passenger volume data (DV and AV, respectively), also obtained from SEA \cite{dvav}, as additional exogenous inputs in our predictive model at each time step. AV and DV for a given time step are determined from the preceding two hours of aggregated incoming volume and the subsequent two hours of aggregated departing volume data, respectively. ~\Cref{fig:pass-volume-norm} depicts the normalized cumulative passenger volume. The resulting input vector $u_k$ is given as, $u_k=\begin{bmatrix} \textrm{TD, TA, DV, AV} \end{bmatrix}^\top$.
To determine the macroscopic traffic dynamics model in~\eqref{eq:dynamic_model}, we train a linear regression model that estimates matrices $A$ and $B$. 

Let $X$ represent a collection of states, $x_k, \ k\in \{1, \ldots, N-1\}$, where $N$ is the size of our dataset. Let $Y$ be the set of states at the following time step,  $x_{k+1}$ and $U$ the set of inputs, whose elements correspond to each entry in $X$. We randomly select 80\% of the data points as the training set, $X_{\mathrm{Tr}}\subset X$ (and $U_{\mathrm{Tr}}\subset U, Y_{\mathrm{Tr}}\subset Y$) to represent the training dataset and reserve the rest for validation. The state and input matrices may be stacked to obtain $X'=[X_{\mathrm{Tr}};U_{\mathrm{Tr}}]$. The resulting system matrix is $A'=[A\ B]\in \mathbb{R}^{4\times 8}$. The optimization problem to estimate the best fit model in~\eqref{eq:dynamic_model} is given by:
\begin{align}\label{eq:Opt1}
    \argmin_{A'}||Y-A'X'||_2 +\rho ||A'||_{2,1},
\end{align}
where $||M||_{2,1}=\sum_{i=1}^n\sqrt{\sum_{j=1}^mm_{ij}^2}$ is the $l_{2,1}$-norm of a matrix first introduced in~\cite{ding2006r} as a rotational invariant $l_1$ norm shown to be less sensitive to outliers~\cite{nie2010efficient}.

We can augment the objective function above with constraints that encode domain knowledge about traffic flow and other system characteristics.
For example, Departures and Arrivals flow are measured on separate and independent roadways and are likely to be uncorrelated. Similarly, treating Arrivals is expected to positively impact the speed at Arrivals (and similarly for Departures). To include this knowledge, we change~\Cref{eq:Opt1} to:
\begin{subequations}\label{eq:Opt2}
\begin{align}
     \argmin_{A'}&||Y-A'X'||_2 +\rho ||A'||_{2,1}\\
     \text{s.t.} \ \ &CA'=0\\
     &DA' \le 0
\end{align}
\end{subequations}
where $C\in \mathbb{R}^{4\times 4}$ encodes the independence of Departures and Arrivals flow, and $D\in \mathbb{R}^{4\times 4}$ encodes the positive impacts of treating roadways on their corresponding speeds. Finally, $\rho$ is a regularization term. The values of these parameters can be found in the code on the \href{https://github.com/pnnl/dynamiccurbs/tree/master/seatac_vms/seatac_vms_controller}{Github page}.

\begin{figure}[t]
    \centering
    \includegraphics[width=\columnwidth]{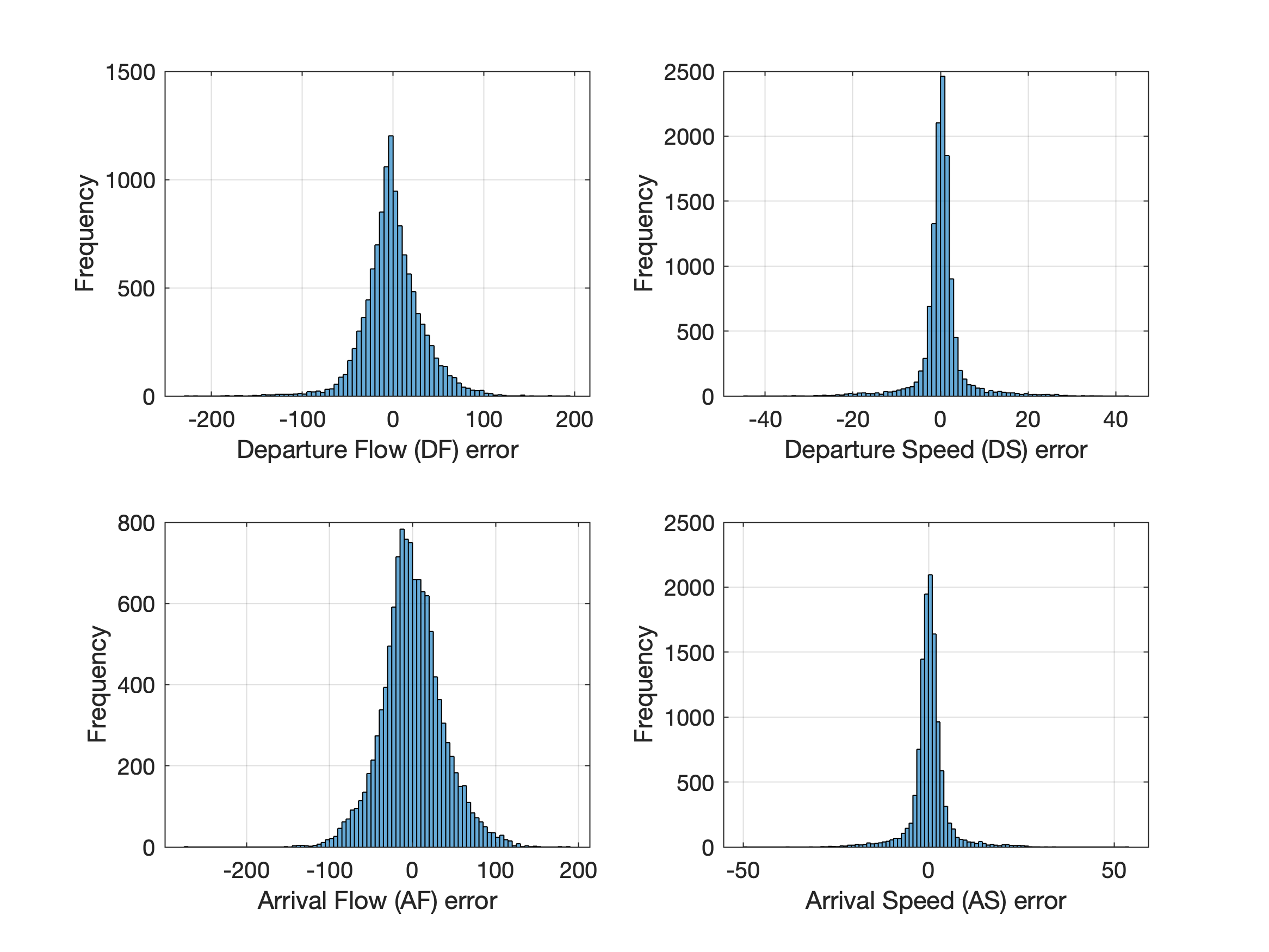}
    \caption{Error histograms for arrival and departure flows and speed using the model obtained from solving~\Cref{eq:Opt2}.}
    \label{fig:model_error_hist}
\end{figure}

The optimization problem in~\eqref{eq:Opt2} is convex and can be solved with a generic optimization solver. We evaluate the prediction accuracy of the resulting $A$ and $B$ matrices against a held-out portion of our ground truth dataset. Each held-out dataset point is a tuple of the previous and next system state. We apply the trained linear dynamics model (consisting of $A$ and $B$) on the previous state to predict the next state, and compute the error between the predicted and true system state at the next time interval. 
This model abets our primary goal of control design for automated decision-making. More accurate models may be built to estimate macroscopic traffic dynamics, but it is beyond the scope of our current research.

\begin{table}[h]
    \centering
    \begin{tabular}{@{} lcc @{}}
    Variable & MAE & RMSE\\
    \midrule
    Departure Flow (count)  &  22.87 & 31.73\\
    Departure Speed (kph)  & 2.70 &  4.99\\
    Arrival flow (count)  &  28.12 &  36.85\\
    Arrival speed (kph) & 2.87 & 4.89\\
    \bottomrule
    \end{tabular}
    \caption{Error statistics on our regression model.}
    \label{tab:error-stats}
\end{table}

~\Cref{fig:model_error_hist} shows a histogram of the errors for each state variable (DF, DS, AF, AS), while~\Cref{tab:error-stats} shows the mean absolute error (MAE) and root mean square error (RMSE) for the same variables. These results indicate that the model is good enough to predict a congested or uncongested scenario and inform the design of TA or TD.
In the following section, we formulate an MPC framework that predicts traffic congestion and designs the timing of the diversionary treatment signals to shorten the time spent in a congested state. We use the distribution of the errors in~\Cref{fig:model_error_hist} to augment the trained model predictions to closely mimic that of the true system and demonstrate the improved performance.

\section{Model Predictive Control for Diversions}
\label{sec:mpc}
We now integrate the macroscopic landside airport traffic model from~\Cref{sec:dynamics} with the following MPC formulation: 
\begin{subequations}\label{eq:MPC}
\begin{align}
    \min_{u_k[1:2]} &\sum_{k=0}^T((\mathcal{D}_k-1)^2+(\mathcal{A}_k-1)^2 + \gamma \sum u_k[1:2])\label{eq:MPC_obj}\\
    \text{s.t} \ \ &\mathcal{D}_k\le 1\label{eq:MPC_C1}\\
    &\mathcal{D}_k \le \frac{x_k[2]}{\widehat{DS}}\label{eq:MPC_C2}\\
    &\mathcal{A}_k\le 1\label{eq:MPC_C3}\\ &\mathcal{A}_k \le \frac{x_k[4]}{\widehat{AS}}\label{eq:MPC_C4}\\
    &x_{k+1}=Ax_k+Bu_k\label{eq:MPC_C5}\\
    &x_0=x_{\text{init}}\label{eq:MPC_C6}\\
    &x_k \ge 0\label{eq:MPC_C7}\\
    &\sum u_k[1:2] \le 1\label{eq:MPC_C8}\\
    &u_k[1:2] \in \{0,1\} \label{eq:MPC_C9}
\end{align}
\end{subequations}
where $\widehat{DS}=35$ km per hour and $\widehat{AS}=45$ km per hour are the critical speeds for Departures and Arrivals respectively, determined from the speed flow fundamental diagrams~\cite{vasisht2022Estimating}. $u_k[1:2]$ denotes the values of TD and TA at time $k$. $x_{\text{init}}$ is the initial state of the system at time $t=0$ and $\gamma$ is the regularization factor for the control signal. $\mathcal{D}_k$ and $\mathcal{A}_k$ are Departures and Arrivals critical speed ratios, defined as the ratio of the average speed to the critical speed. \textit{A critical ratio value well below $1$ indicates significant congestion}.

The objective \eqref{eq:MPC_obj} is to minimize the deviation of $\mathcal{D}_k$ and $\mathcal{A}_k$ from $1$. The constraints in~\eqref{eq:MPC_C1}-\eqref{eq:MPC_C4} penalize Departures and Arrivals speeds only when the respective critical ratios fall below $1$, but not when they are above $1$, which is consistent with our goal of increasing average speeds (and hence reducing congestion). Additional constraints, such as system dynamics, initial state and restricting the state variables to be positive are described in ~\eqref{eq:MPC_C5} -~\eqref{eq:MPC_C7}. Finally, the constraints in~\eqref{eq:MPC_C8}-\eqref{eq:MPC_C9} limit the control inputs (TD and TA) to be binary variables and that at most only one diversionary signal is active at each time step (since the VMS can only display one message at a time).

We solve the problem in~\Cref{eq:MPC} in a receding horizon manner, i.e., at each time step, the above objective is minimized over an appropriate time horizon to obtain the optimal sequence of TD or TA signals. The first elements of this sequence are applied as input to the system and the dynamics is propagated forward in time. This process is repeated by resetting $x_0$ with $x_{k+1}$ until the final time step. Since we do not have access to real-time measurements from the actual system or a simulator to represent it, we use the distributions of the model errors shown in~\Cref{fig:model_error_hist} to generate system outputs as a response to control (or no control). The optimal control signals obtained from solving~\Cref{eq:MPC} are applied to the model:
\begin{align}\label{eq:random_model}
    x_{k+1}=Ax_k+Bu_k+\mathcal{U}(a,b)
\end{align}
where $\mathcal{U}(a,b)$ is a uniform distribution between $a$ and $b$, based on the 90th percentile values of the errors shown in~\Cref{fig:model_error_hist}. The values of $a$ and $b$ for the four states are provided in Table~\ref{tab:uniform-dist}.

\begin{table}[b]
    \centering
    \begin{tabular}{@{} lcc @{}}
    Variable & a & b\\
    \midrule
    Departure Flow (DF) error  &  -34.6 & 38.5\\
    Departure Speed (DS) error  & -2.9 &  3.4\\
    Arrival flow (AF) error  &  -43.3 &  46.5\\
    Arrival speed (AS) error & -3.5 & 3.9\\
    \bottomrule
    \end{tabular}
    \caption{Parameters $a$ and $b$ describing the uniform distribution $\mathcal{U}(a,b)$ from~\Cref{eq:random_model}.}
    \label{tab:uniform-dist}
\end{table}

We demonstrate the effectiveness of our formulation by performing several runs of the optimization loop with the randomized system output from~\Cref{eq:random_model} and recording the mean and standard deviation of the resulting speeds and flows. We now discuss how we use these statistics to compute metrics such as travel time savings and speed improvement.

\section{Evaluation}
\label{sec:eval}

\begin{figure}[t]
    \centering
    \includegraphics[width=\columnwidth]{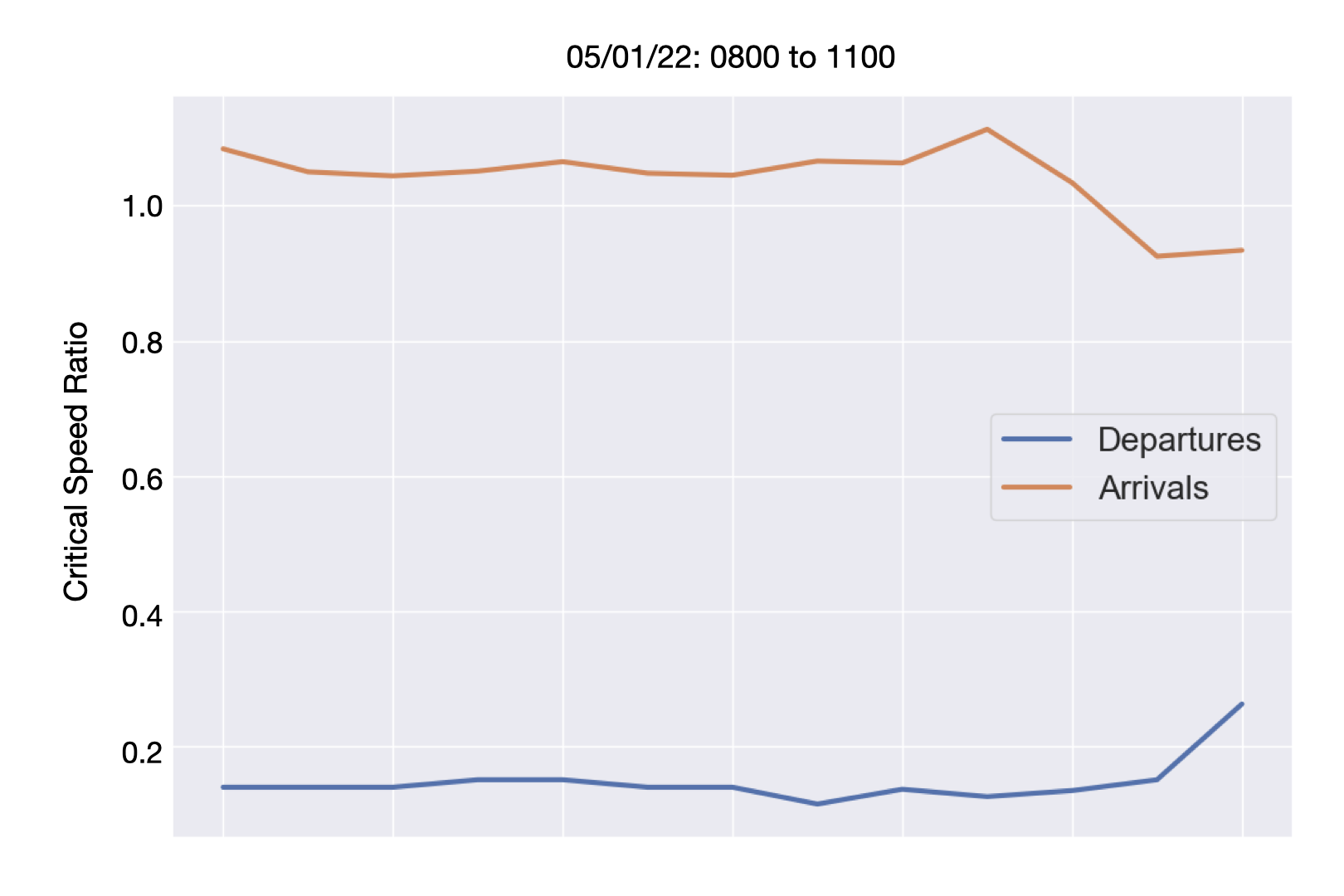}
    \caption{An untreated congestion scenario from our dataset. On May 1, 2022, the Departures roadway (blue line) at SEA was congested for nearly three hours from 0800-1100, while Arrivals was operating normally.}
    \label{fig:untreated_cong}
\end{figure}

This section describes how we evaluate the effectiveness of our approach in recommending diversions to ease congestion. We selected $50$ historical scenarios from our SEA dataset that represent \emph{untreated congestion}, i.e., where either Arrivals or Departures was congested but no diversion was deployed (the scenarios were evenly split between those with Arrivals congested and those with Departures congested). To extract such scenarios, we looked for periods of time, typically 2-3 hours, where the critical speed ratio of either facility dropped significantly while that of the other facility was normal, e.g.,~\Cref{fig:untreated_cong}. These $50$ scenarios serve as the grounding examples on which we run our controller and measure how much congestion would have been eased.

We run our controller on each untreated congestion scenario to yield a sequence of control actions that encode either diversionary messages or no intervention. For any MPC approach, two key parameters are the planning horizon (inner loop) and the number of execution steps (outer loop). Our planning horizon is $12$ time-steps, or $3$ hours, as congestion typically builds up and dissipates within that time-frame. However, we set the number of execution steps to $4$, or $1$ hour. For each scenario, we diverge from history after the first execution step,  since we use our one-step dynamics model to simulate the system state for the next execution step. Therefore, we limit execution steps to $4$ to allow enough time-steps to evaluate our controller meaningfully, but not so much as to diverge significantly.

\begin{figure*}[t]
    \centering
    \begin{subfigure}{\columnwidth}
        \centering
        \includegraphics[width=\columnwidth]{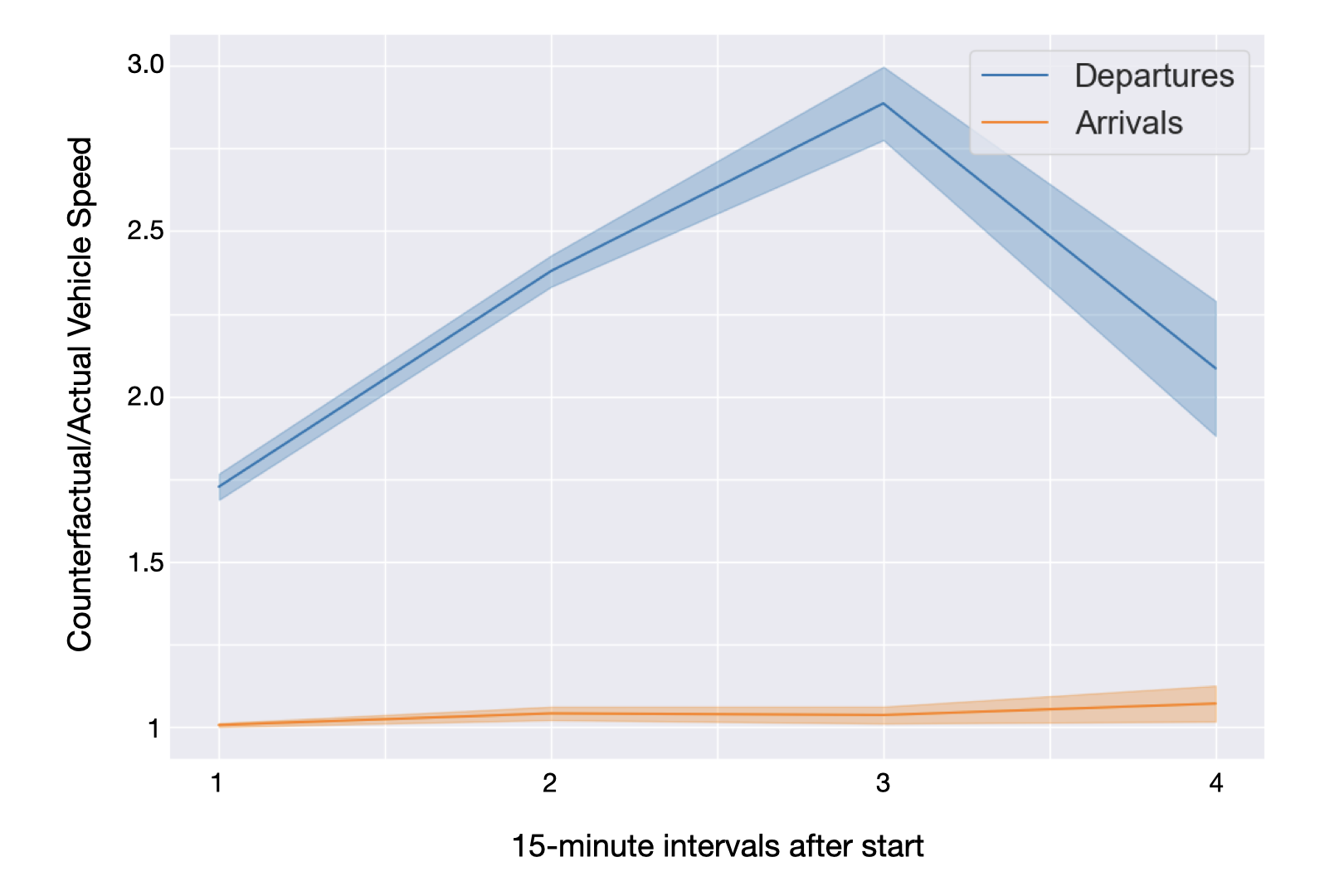}
        \caption{Departures Treated}
        \label{fig:eval-speed-dep}
    \end{subfigure}
    \begin{subfigure}{\columnwidth}
        \centering
        \includegraphics[width=\columnwidth]{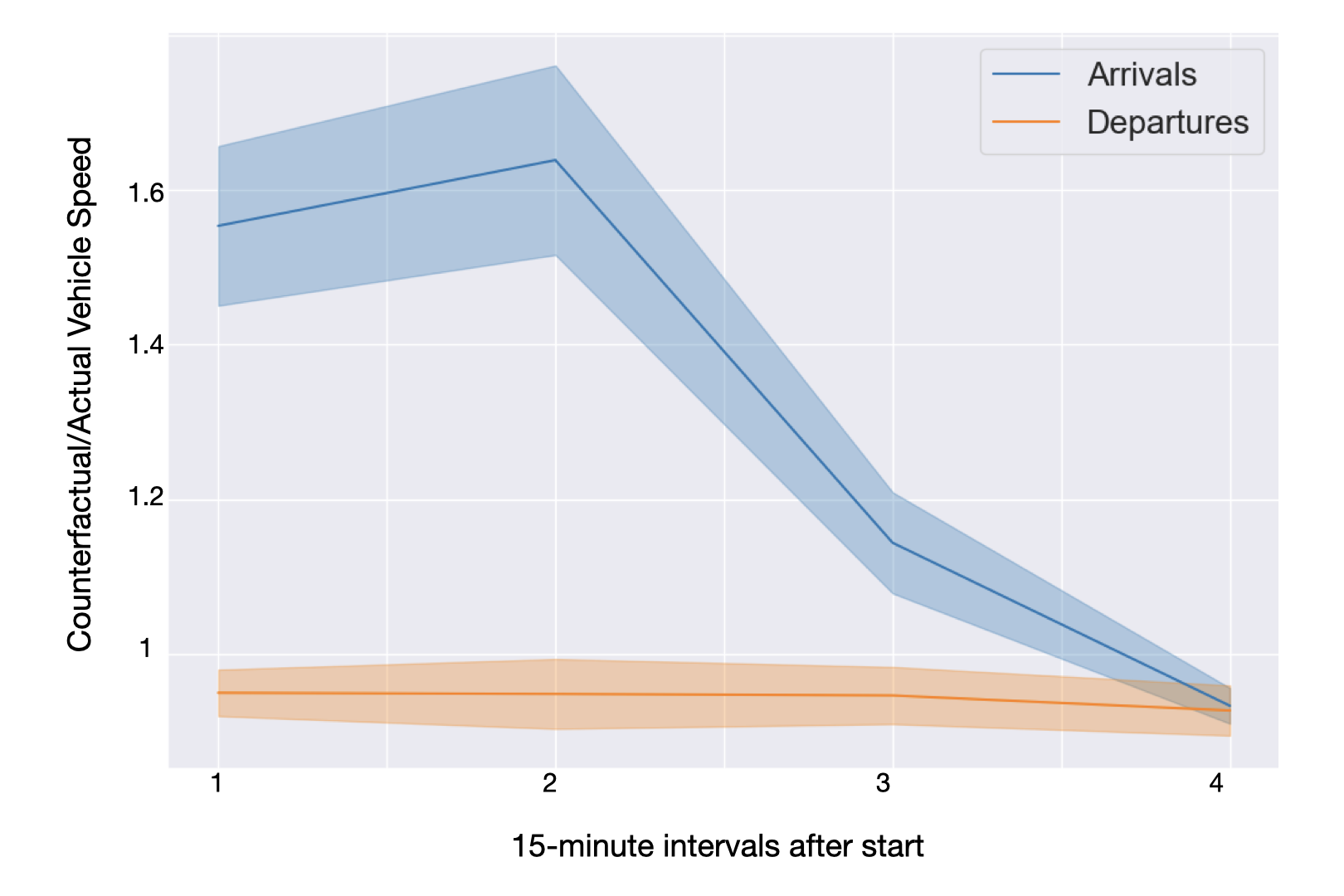}
        \caption{Arrivals Treated}
        \label{fig:eval-speed-arr}
    \end{subfigure}
    \caption{Our MPC controller computes diversions that, if followed, would have improved speed in the treated roadway up to three times the actual speed in the historical scenario. The speed improvement for treating Departures is much higher because our dynamics model estimates a higher impact of treating Departures on its speed, than for treating Arrivals.}
    \label{fig:eval-speed}
\end{figure*}

In addition to a sequence of $4$ control actions, the controller also outputs a sequence of $4$ \textit{counterfactual states} (flow, speed etc.), for each scenario. These values represent the expected future states of the system, up to $1$ hour, had our controller's diversions been deployed. Therefore, we use the actual and counterfactual states of each scenario to obtain two key performance metrics.

The first performance metric uses the actual and counterfactual vehicle speeds. Recall from~\Cref{eq:MPC_obj} that our optimization objective is to improve the speed of the treated roadway without significantly worsening the other roadway. Therefore, we compute the ratios of the counterfactual and actual speeds for each roadway,  up to $4$ execution steps, averaged over the $50$ scenarios. \Cref{fig:eval-speed} \textit{shows that our diversions would have improved speed in the treated roadway by up to three times} without affecting the other roadway. 

The speed improvement ratios for treating Departures (\Cref{fig:eval-speed-dep}) are considerably higher than for treating Arrivals (\Cref{fig:eval-speed-arr}). The primary reason for this difference is our dynamics model from~\Cref{sec:dynamics}, which we use to simulate the counterfactual states. The model estimates a higher impact of treating Departures (TD) on improving Departures speed, than the impact of treating Arrivals (TA) on improving Arrivals speed. The higher impact of TD is in turn due to the higher response rates of drivers to the corresponding diversionary message; our prior work on inferring the response rates explores their differences in detail~\cite{vasisht2022Estimating}.

Note also from~\Cref{fig:eval-speed-arr} that when Arrivals is treated, the standard error around the average counterfactual/actual speed ratio reduces with increasing 15-minute intervals. This reduction contradicts both our intuition that the counterfactual outcomes should become noisier as we simulate further into the future and the fact that this intuition is reflected by the Departures line in~\Cref{fig:eval-speed-dep}. The reason is that congestion in Arrivals typically dissipates more quickly than for Departures, usually within an hour from when congestion is highest. Therefore, in the untreated congestion scenarios used for~\Cref{fig:eval-speed-arr}, the actual speed goes back to hovering around or above the critical speed within the hour. Our controller's diversions seek to bring the counterfactual speed to the critical speed, thus explaining why the ratio of counterfactual to actual speeds is consistently close to 1 over all such scenarios.

\begin{figure}
    \centering
    \includegraphics[trim={1cm, 1cm, 1cm, 1cm},clip,width=\columnwidth]{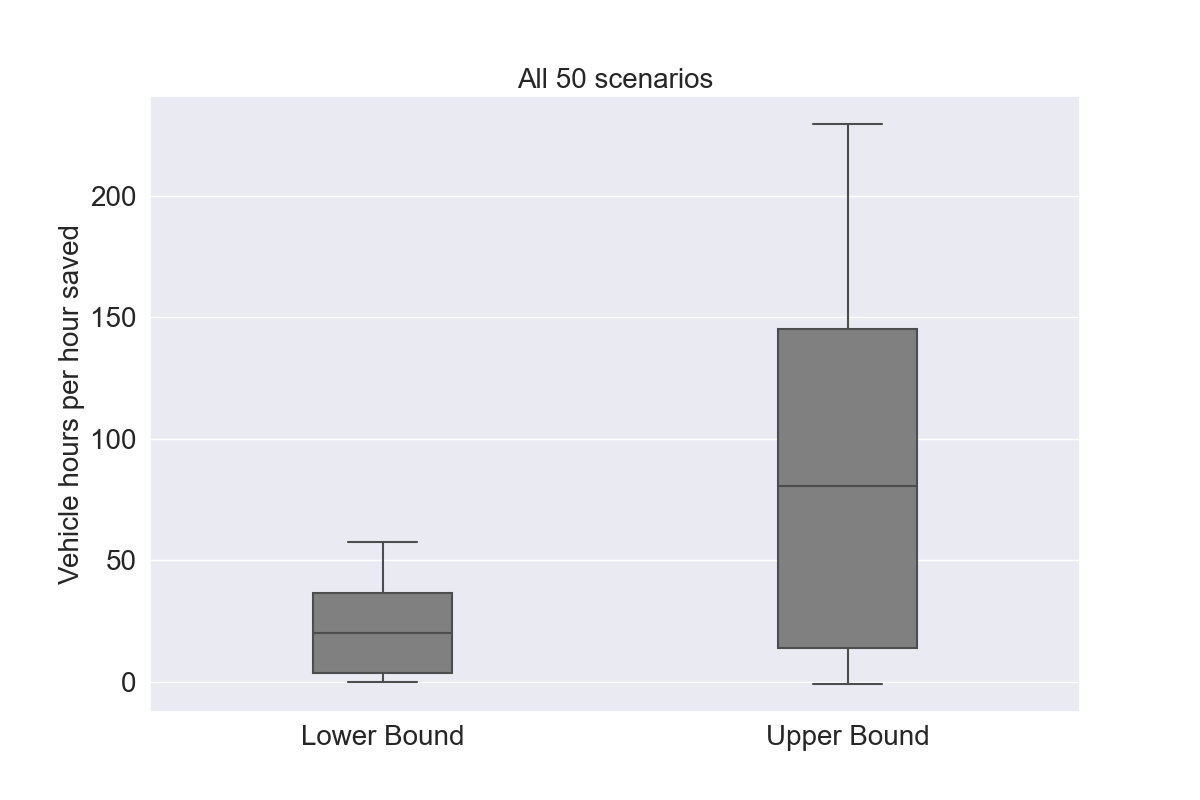}
    \caption{The distributions of the vehicle hours saved due to our controller's diversions, over the $50$ evaluation scenarios. The upper bound for each scenario is $4$ times the corresponding lower bound, given the range of lengths of plausible traffic bottlenecks at Seattle-Tacoma airport.}
    \label{fig:eval-hrs}
\end{figure}

The second performance metric uses both speed and flow; it computes the counterfactual cumulative vehicle hours saved by the vehicles that go to the treated roadway. At SEA, the traffic bottleneck can be anywhere from $0.5$ to $2$ kilometres in length, yielding corresponding lower and upper bounds of travel time improvement respectively. We compute the dot product of the travel time improvement in each execution step with the number of incoming vehicles in the corresponding $15$ minute interval, to get the metric of cumulative vehicle hours saved (per hour).

\Cref{fig:eval-hrs} plots the means and inter-quartile ranges of the lower and upper bounds of the cumulative vehicle hours saved, over the $50$ scenarios. It shows that \textit{our diversions could have saved roughly $20$ to $80$ vehicle hours per hour of deployment} on average. Reducing vehicle hours also reduces idle time fuel wastage and emissions, both important downstream operational metrics. The specific savings depend on vehicle make and model, but we can derive coarse estimates for a mid-sized sedan~\cite{gaines2013greener}. The average idle time fuel expense is $0.35$ gallons per hour per vehicle, which translates to a savings of $7$ to $28$ gallons of fuel per hour. The average idle time emissions is $2100$ grams of carbon dioxide per vehicle per hour, which translates to a total savings of $90$ to $360$ kilograms of $\text{CO}_2$ per hour. Given that the airport experiences around $6$ to $8$ hours of congestion every business day, all these savings could accumulate significantly over business months and years.

\section{Conclusion}
We presented a comprehensive data-driven framework for controlling variable message signs to reduce ground traffic congestion and thus travel time at airports. Using available real-world speed and flow data from SEA, we obtain a linear dynamical model of the system that we then use in a model predictive controller to optimize for our congestion objective. Our simulation results on several historical congested scenarios at SEA showcase the efficacy of our approach and how algorithmic decision-making can augment human judgment in real-world transportation problems. Our specific application requires only one control decision per time-step, but our methodology is general and capable of solving more complex and multi-dimensional control problems.

Future efforts will work towards a pilot deployment of this controller at SEA, to automate the traffic diversion messages and quantify the benefits to travelers in terms of travel time saved. Future work will also explore the human response to such variable message signs and how that can be included within the predictive control framework.

\section*{Acknowledgements}
We thank the landside operations team at Seattle-Tacoma International Airport and Port of Seattle for sharing data and insights.
\label{sec:concl}

\printbibliography

\end{document}